\documentclass[12pt]{article}
\usepackage{graphicx}
\newcommand{\EE}{{\rm I\kern-2pt E}}
\newcommand{\RR}{{\rm I\kern-2pt R}}
\newcommand{\DD}{{\rm I\kern-2pt D}}
\newcommand{\PP}{{\rm I\kern-2pt P}}
\newcommand{\NN}{{\rm I\kern-2pt N}}

\title{THE SHIFT:\\ Properties and recommendations for practical use}
\author{Nicolas BOULEAU \\ {\it Ecole des ponts, ParisTech}} \date{}
\begin{document}
\maketitle

\noindent {\bf Abstract.} In mechanics of structures as well as in other
domains of engineering, the probabilistic models often have to be computed by
simulation. One of the reasons is that the stochastic calculus involved in non
linear representations yields rarely explicit formulas. Usually these
simulations use, for each sample, a large number of calls to the random
function, in this case the simulation by the shift, whose field
of application is as wide as that of  Monte Carlo method, is particularly
relevant.

The theoretical features, the implementation and the specific advantages of
this method have been taught in a course of the author at Paris VI University
in 1988 and are detailed in the book with D. L\'{e}pingle [BL].\\

\noindent{\bf I. Other known methods of simulation}

Let us begin by the presentation of the methods to be compared with the shift.\\

\noindent{\bf A. The Monte Carlo method}

For computing the expectation $\EE(X)$ of a random variable $X$, it consists of

a) representing $X$ as a random variable defined on the probability space
$$([0,1]^s, {\cal B}[0,1]^s, dx_1dx_2\cdots dx_s)$$ with $s$ finite or infinite,

b) imitating independent random samples of points of $[0,1]^s$,
$$(U_{11},U_{12},\ldots,U_{1s})\;,\;(U_{21},U_{22},\ldots,U_{2s})\;,\ldots,\;(U_{n1},U_{n2},\ldots,U_{ns})\;,\ldots$$

c) and applying the law of large numbers
$$\EE(X)=\lim_{N\uparrow \infty}\frac{1}{N}\sum_{n=1}^N
X(U_{n1},U_{n2},\ldots,U_{ns}).$$

\noindent{\bf Comments}

\noindent$\bullet$ The step a) is the art of using the {\it algorithms of
simulation}. It is theoretically possible to take $s=1$, but practically, the
dimension $s$, finite or infinite, is rather naturally given and is difficult
to reduce. The main reference for these algorithms is [Devroye 86]. It is worth
to note that  this operation of representing a random variable $X$ on the cube
$[0,1]^s$ ($s$ finite or not) can be effectively performed in very much
more cases than those where the law of $X$ is known by an explicit formula.
This is, in particular, due to the {\it rejection method} which allows an exact
simulation of a random variable even when its density is only known by a
sequence of approximations.

\noindent$\bullet$ The step b) is the generation of pseudo random numbers. The
imitation of randomness cannot be perfect. It has been theoretically proved by
logicists, especially by Martin L\"{o}f. Practically extremely good generators
are available with gigantic period. They are obtained by testing the production
algorithms through statistical tests and eliminating those who behave badly.
The tests are quite numerous (see [Niederreiter 1978], Knuth 1981, Marsaglia
1985, Ripley 1987, Fushimi 1988, Altman 1988, Anderson 1990, and other
references given in [BL])

\noindent$\bullet$ The class of functions to which the Monte Carlo method
applies is theoretically all functions in ${\cal L}^1([0,1]^s,dx)$ (here $dx$
is the product probability measure). Practically, it is safe to restrict the
method to bounded functions. Nevertheless the irregularity of the function is
not limited.

\noindent$\bullet$ The speed of convergence is given for $X$ in $L^2$, hence
for $X$ bounded, by the law of iterated logarithm
$$\limsup_N \frac{1}{\sigma\sqrt{2N\log\log N}}\left(\sum_{n=1}^N
X(U_{n1},U_{n2},\ldots,U_{ns})-\EE X\right)=1$$ where $\sigma^2$ is the variance
of $X$.

There are also global estimates of the type Cramer-Chernov ($X$ bounded)
$$\PP(\frac{1}{N}\sum_{n=1}^N X(U_{n1},U_{n2},\ldots,U_{ns})\geq \EE
X+\varepsilon)\leq e^{n\psi(\varepsilon)}$$ which are in connexion with large
deviation theory.\\

\noindent{\bf B. Quasi-Monte Carlo methods}

The step a) is preserved, but b) and c) are replaced by

b') choosing an equidistributed sequence $(\xi_n)$ on $[0,1]^s,$

c') putting $$\EE(X)=\lim_{N\uparrow\infty}\frac{1}{N}\sum_{n=1}^N X(\xi_n).$$

\noindent{\bf Comments}
\noindent$\bullet$ This method requires the random variable
to be represented on $[0,1]^s$ by a {\it Riemann
integrable function} i.e. a bounded function whose
discontinuity points belong to a negligible set, or
equivalently such that $\forall \varepsilon>0$ continuous
functions $u$, $v$ exist such that $u\leq X\leq v,
\;\int(v-u)dx\leq \varepsilon$.

Often the function obtained on $[0,1]^s$ by the step
a) is Riemann integrable already, so this fact brings no
restriction.

The non-Riemann integrable functions encountered in
practice, come from stochastic calculus related to
Brownian motion for which quasi-Monte Carlo methods are
not yet relevant.

\noindent$\bullet$ The speed of these methods using low
discrepancy sequences [cf Niederreiter 1978] is impressive in low
dimension. A great advantage is also that the rate of convergence is
given by a deterministic criterion --- the Koksma-Hlawka formula ---
allowing to produce the numerical value of the computed expectation
with an explicit accuracy (cf [BL] Chapter 2C). But this speed is
decreasing when the dimension increases, so that for large dimensional
simulations these methods are at present irrelevant (see [BL] for the
study of the realm of efficiency of these methods).\\

\noindent{\bf II. The shift method}\\

\noindent{\bf A. The infinite dimension in practice}

Because the shift method is particularly relevant in large or infinite
dimension, it is important to emphasize the fact that {\it the case of
infinite dimension is practically the most frequent}.

Indeed, most often, in the simulation, the number of calls to the
random function, although finite, is itself random and unbounded. A
typical example is the computation  of the expectation of a stopping
time of a Markov chain, e.g. the entrance time in a given set.

But even for a finite dimensional probabilistic model the step a) will
lead to the infinite dimension if the rejection method is used, and
practionners know how rejection is often unavoidable.\\

\noindent{\bf B. The principle of the shift method}

It is based on the pointwise ergodic theorem instead of the law of
large numbers.

Let $X$ be a random variable simulated under the form
$$X=F(U_1,U_2,\ldots,U_n,\ldots)$$ where yhe $U_i$'s are i.i.d. random
variables uniformly distributed on $[0,1]$. In other words the $U_i$'s
are the coordinate mappings from $[0,1]^\infty$ equipped with the
product Lebesgue measure to each factor. The method consists of putting
$$\EE(X)=\lim_{N\uparrow\infty}\frac{1}{N}[F(U_1,U_2,\ldots)+F(U_2,U_3,\ldots)+\cdots+F(U_N,U_{N+1},\ldots)]$$
instead of using in the case of the Monte Carlo method a double
sequence of calls:
$$\EE(X)=\lim_{N\uparrow\infty}\frac{1}{N}[F(U_{11},U_{12},\ldots)+F(U_{21},U_{22},\ldots)+\cdots+F(U_{N1},U_{N2},\ldots)]$$
The theorem of Birkhoff says that the shift method converges almost
surely and in $L^1$ as soon as $X$ is in $L^1$ and in $L^p$ as soon as
$X$ is in $L^p$ $1\leq p<\infty$.\\

\noindent{\bf C. General implementation manner}

Suppose we have written a procedure able to give us a sample of the
probabilistic model we are studying. 

Instead of writing the main program using this procedure successively
$N$ times and doing the average, we shall write the main program using
$N$ times the procedure after initialisation and shifting one step the
random function each time.

The economy on the pseudo-random numbers generator is evident.\\

\noindent{\bf D. Fine implementation for efficiency}

In order to exploit the full strength of the method the maximum amount
of information of the preceding sample has to be kept for the following
one.

This can be done very efficiently by {\it pointers}. Let us take the
example of a Markov chain for the explanation.

Let us consider a Markov chain $X_n$ with values in $\RR^d$ which is
simulated on $[0,1]^\infty$ in the following way:
$$X_{n+1}=F(X_n,n,U_{n+1}),\quad X_0=x$$where $F$ is a map from
$\RR^d\times \NN\times[0,1]$ into $\RR ^d$, and where the $U_n$'s are,
as before, the coordinate mappings of the cube $[0,1]^\infty$ to its
factors.

Let us suppose we want to compute the expectation of a functional of
the process $X$, for instance $\EE[G(X_T,T)]$ where $G$ is a given
bounded function and $T$ the hitting time in the set $A\subset\RR^d$:
$T=\inf\{n>0;\;X_n\in A\}$.

Denoting $\theta$ the shift operator on $[0,1]^\infty$ defined by
$$U_n\circ\theta=U_{n+1},$$ the ergodic theorem writes
$$\EE[G(X_T,T)]=\lim_{N\uparrow\infty}\frac{1}{N}\sum_{n=0}^{N-1}G(X_T,T)\circ\theta^n.$$
The use of pointers can be explained as follows. In order to apply the
preceding formula we have to compute$G(X_T,T)$ on successive points
$\omega,\theta(\omega),\theta^2(\omega),\ldots,$ in $[0,1]^\infty$.
These points are sequences of points in $[0,1]:$ $$\begin{array}{rcl}
\omega &=&(U_1(\omega),U_2(\omega),\ldots,U_k(\omega),\ldots)\\
\theta(\omega) &=&
(U_2(\omega),U_3(\omega),\ldots,U_{k+1}(\omega),\ldots)
\end{array}$$
If, as we suppose, the stopping time $T$ is finite, we need only to
compute a finite length of each of these sequences: the sequence
$\omega$ is computed until $T(\omega)$ which is given by the test
$$<\mbox{test}>\left\{\begin{array}{rcl}
X_k(\omega)\mbox{ still outside }A & \Rightarrow & k<T(\omega)\\
X_k(\omega)\in A\mbox{ for the first time }&\Rightarrow &k=T(\omega)
\end{array}\right.$$
Suppose we have picked out the sequence $\omega=(U_1,U_2,\ldots)$ and
we have put the numbers $U_1,U_2,\ldots$ in pointers as in the following
figure:

\begin{center}
\includegraphics[width=4in]{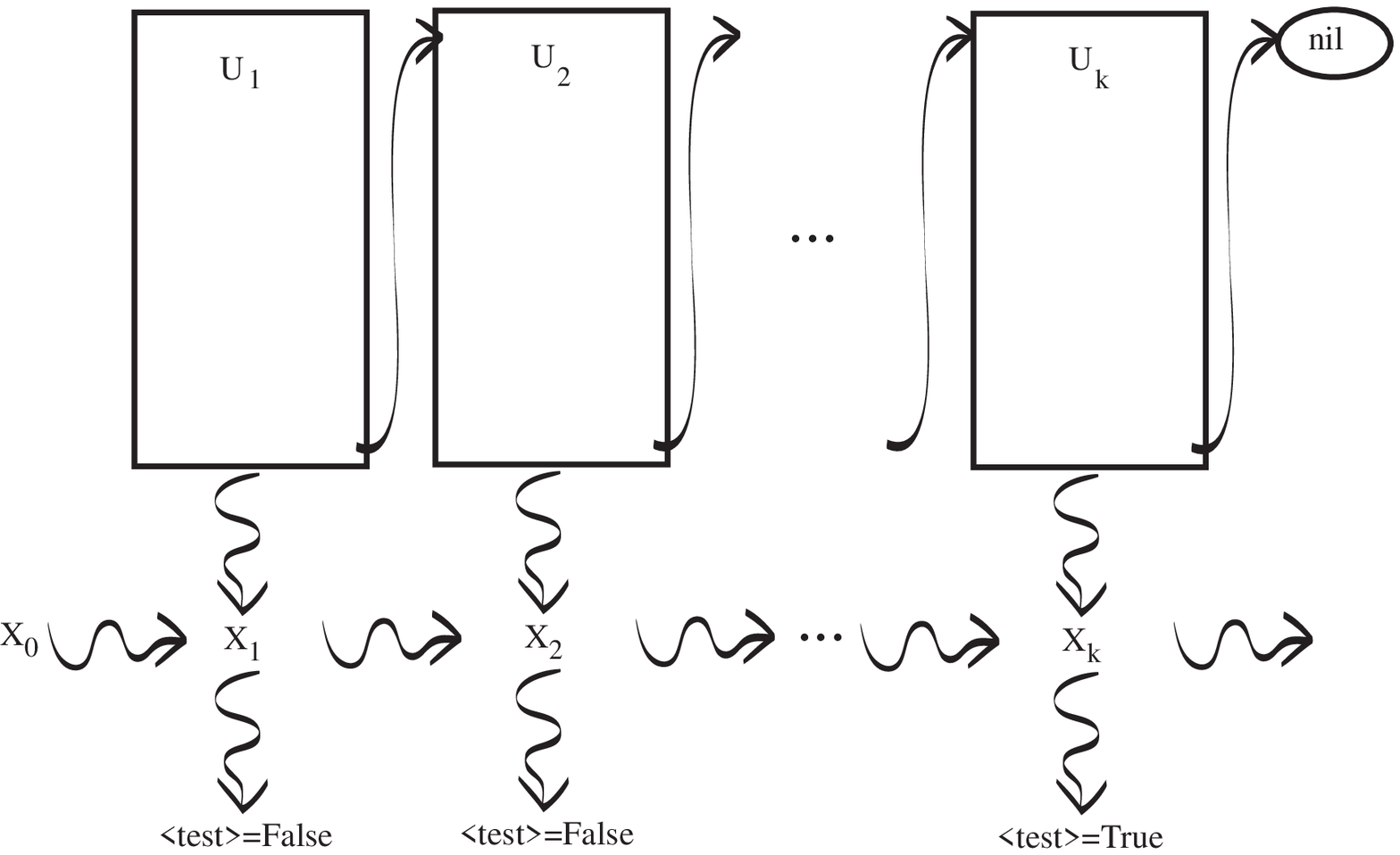}
\end{center}

\noindent The sequence $\theta(\omega)=(U_2,U_3,\ldots)$ is already partially
chosen: either it is long enough or it must be lengthened, in which case
we have the new scheme:

\begin{center}
\includegraphics[width=5in]{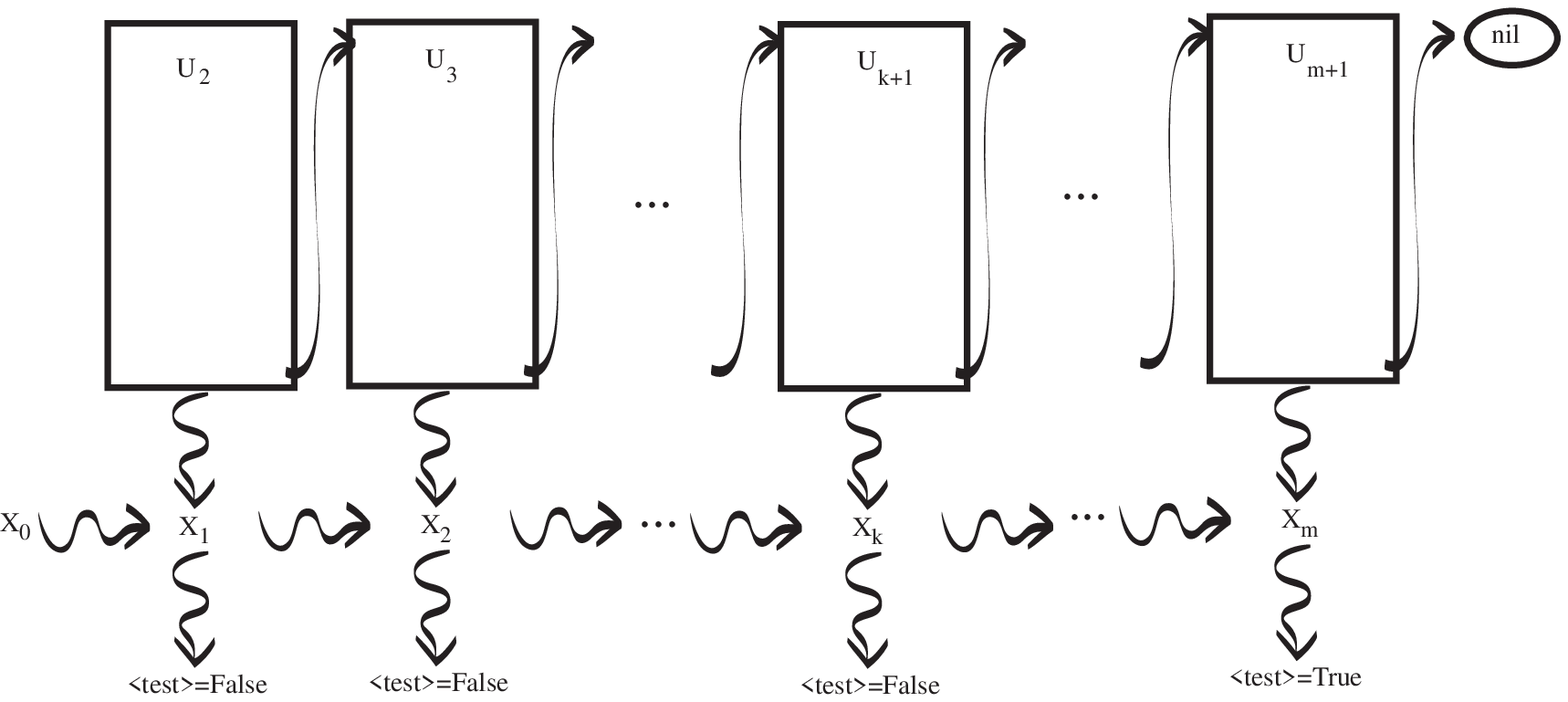}
\end{center}

\noindent In these figures the arrows represent computations to do,
essentially to compute
$$X_{n+1}=F(X_n,n,U_{n+1})$$ When we compute $F(x,n,y)$ however, only
the variable $x$ is new (and $n$ becomes $n+1$) and {\it partial
computations depending only on $y$ can be stored in the box of $U_n$}.
This storage must of course be made while we are lengthening the
sequence as explained above.

Clearly this storage of partial computations cannot be done in the
Monte Carlo method.\\

\noindent{\bf E. Theoretical results on the rate of convergence}

We dont go into the details of these results here, see  [BL] for mathematical
proofs. 

There is no standard rate of convergence valid  for every function in $L^2$.
Nevertheless in practice usual functions belong to a class called the Gordin
class in which a law of iterated logarithm holds.

This law of iterated logarithm involves a coefficient different from the
variance, which can be either bigger or smaller than the variance even in the
finite dimensional case. So even for integration in finite dimension the shift
can be faster than Monte Carlo.

From a practical point of view the possibility of storage of partial
computations is the dominant phenomenon.\\

\noindent{\bf F. Numerical example}

The simulated device is a transport particles problem. A particle arrives in
$0$ coming from the real negative axis and enters into the square
$[0,1]\times[-\frac{1}{2},\frac{1}{2}]$ where it goes in straight lines during
a random distance which follows an exponential law with parameter $\lambda$,
then it splits into two similar directions uniformly distributed on
$[-\frac{\pi}{2},\frac{\pi}{2}]$ with respect to the direction of the preceding
particle, and they behave as the initial particle provided they are inside the
square. The quantities to be computed are the mean number of splittings and the
mean number of particles leaving the cube by the righthand side.

The following results have been obtained by simulating 500 000 samples by the
Monte Carlo method (MC) and the shift method (Sh)\\

{\small\begin{tabular}{r||lr||lr||lr||lr||lr||}
parameter & \multicolumn{2}{c||}{0.980}
& \multicolumn{2}{c||}{0.960}
& \multicolumn{2}{c||}{0.940}
& \multicolumn{2}{c||}{0.920}
& \multicolumn{2}{c||}{0.900}\\
\hline
method  
& MC & Sh & MC & Sh & MC & Sh & MC & Sh & MC & Sh\\
\hline
mean number&&&&&&&&&&\\
of splittings&2.78&2.79&2.92&2.93&3.12&3.11&3.34&3.31&3.60&3.58\\
\hline
mean number of&&&&&&&&&&\\
particles through&3.97&4.00&4.36&4.38&4.95&4.90&5.68&5.62&6.72&6.65\\
the righthand side&&&&&&&&&&\\
\hline
mean number of&&&&&&&&&&\\
calls to the&24.5&4.9&27.4&4.9&32&4.9&37&5&46&5\\
random function&&&&&&&&&&\\
\hline
ratio of the&&&&&&&&&&\\
duration of the
&\multicolumn{2}{c||}{2.04}
&\multicolumn{2}{c||}{2.19}
&\multicolumn{2}{c||}{2.31}
&\multicolumn{2}{c||}{2.40}
&\multicolumn{2}{c||}{2.46}\\
programs MC/Sh&&&&&&&&&&\\
\hline
\end{tabular}}\\

\vspace{.5cm}
\noindent{\bf III. Recommendations}

To become an expert in the art of using the shift, it is essential to keep in
mind that {\it the rate of convergence depends on the order of the calls}

Let us take a finite dimensional example: Let $f$ be a bounded function from
$[0,1]^d$ into $\RR$. The integral of $f$ can be obtained by the shift by
putting
$$\EE
f=\lim_{N\uparrow\infty}\frac{1}{N}[f(U_1,\ldots,U_d)+f(U_2,\ldots,U_{d+1})+\ldots+(U_N,\ldots,U_{N+d})]$$
Now we can also permute the coordinates defining with the permutation $\sigma$
$$f_\sigma(x_1,\ldots,x_d)=f(x_{\sigma(1)},\ldots,x_{\sigma(d)})$$ and apply the
shift to $f_\sigma$.

The rate of convergence {\it will be different} in general, so that there are
$d!$ different manners of applying the shift.

In finite dimension these different manners give roughly speaking the same order
of rate (see [BL] for a detailed study)

But in infinite dimension it is quite different and for simulating the
expectation of a functional of a random process this becomes very important:
{\it The idea of discretization of the process and shifting along the time is
one of the {\bf worse} way}. Good implementations for functionals of Brownian motion
and for solution to stochastic differential equations are given in [BL Chapter
V].\\

\noindent REFERENCES

\noindent[BL] Bouleau N. and L\'{e}pingle D.(1994)

 {\it Numerical methods for stochastic
processes}, J. Wiley and Sons

\noindent Devroye, L. (1986)

{\it Non-uniform random variate generation} Springer

\noindent Niederreiter H. (1978)

	Quasi-Monte Carlo methods and
pseudo-random numbers,	

{\it Bull. Amer. Math. Soc.} 84, 957-1041 	 

\noindent Bouleau N. (1991)

	On numerical integration by the shift and
application to Wiener space, 	

{\it Acta Applic. Math.} 25, 201-220

\end{document}